\documentclass[twoside,12pt]{amsart}
\usepackage{myamsart}
\input{mydef}

\begin{document}
\title[No More Than Mechanics. I]{No More Than Mechanics. I%
\\ \large\itshape Plain Mechanics:\\
Classical and Quantum Mechanics as well}
\thanks{This research was partially supported by grant of {\em
Foundation PRO MATHEMATICA\/} (French Mathematical Society).
Final preparation of this work was partially supported by grant of the
CONACYT Project 1821-E9211, Mexico.}
\thanks{On leave from the Odessa State
University.}

\author{Vladimir V. Kisil}

\email{kisilv@maths.leeds.ac.uk}
\address{%
School of Mathematics,
University of Leeds,
Leeds LS2\,9JT,
UK}

\begin{abstract}
This is the written version of a short talk on
{ \em $10^{th}$ Conference on Problems and Methods
in Mathematical Physics\/ } (September 13 - 17, 1993 in Chemnitz,
Germany).
A new scheme of the quantization is presented.
A realization of the scheme for a particle in $n$-dimensional space
by two-sided convolutions on the  Heisenberg group is constructed.
\AMSMSC{81P05, 81S99}{22E25}
\end{abstract}
\maketitle

\tableofcontents

\section{Introduction}

The main goal of this paper is to present a new approach to
relationships between classical and quantum mechanics. I think, there is
a large discrepancy between the mathematical beauty of quantum mechanics
and its absurdity from the common point of view. This discrepancy does
not prevent us to make our computations with a large preciosity, but it
induces  future investigations of our basic assumptions.
The problem of a quantization is still under the serious investigation
(see
\cite{Berezin74,Berezin75,BoutetGuill85,Coburn90,Guillemin84,GuillStern81,Unterberger84})
beyond more than a half century after the creation of
quantum mechanics.

This paper makes principal steps in such direction (but does not
achieve the final point, of course). The usual ``quantization'' means
some (more or less complete) set of rules for the construction of a
quantum algebra from the classical description of a physical system. Our
main suggestion is to replace such quantization by the constitution of
an operator algebra from which both the classical and the ``usual'' quantum
descriptions may be derived. The paper gives the more precise
formulation for this approach and illustrates it on the simplest example: the
quantization for a particle in $n$-dimensional space. Future papers in
this series will present a more concluded description of {\em Plain
Mechanics}. For example, application of plain mechanics to the quantum
field theory requires consideration of Clifford valued convolutions on
the Heisenberg group (see~\cite{Kisil93c}).

I am glad to express my gratitude to the following peoples for their helpful
discussion: Dr.~M.~Kuzmin, Prof.~Yu.~Gurevich, Prof.~N.~Vasilevski and
Prof.~B.~Veytsman.

\section{Mathematical Background}\label{se:background}

Our scheme is based on the properties of convolutions on the Heisenberg
group\footnote{More general case of arbitrary step 2 nilpotent Lie
group may be considered also, but we do not touch this theme here. For
corresponding results see
\cite{Kisil92,VasTru92}.}.
 This subject is well known and there are many suitable sources on it
\cite{MTaylor84,MTaylor86}. So we can introduce a few definitions only.
 The Heisenberg group is a step 2 nilpotent Lie group. As a
$C^\infty -$manifold it coincides with $\Space{ R}{2n+1}$. If an
element of it is given in the form $g=(u,v)\in \Space{H }n$, where
 $u\in \Space{ R}{}$ and $v=(v_1,\dots, v_n)\in \Space{ C}n$, then the
group law on $
\Space{ H}{n}$ can be written as
\begin{displaymath}
(u,v)*(u',v')=\left(u+u'-2 \object{Im} \sum_1^n v'_k \bar{v}_k,\
v_1+v'_1,\ldots, v_n+v'_n\right) .
\end{displaymath}
We single out on $\Space{ H}n$ the group of nonisotropic dilations
 $\{\delta_\tau\}$, $\tau\in \Space{ R}{}_+$:
\begin{displaymath}
\delta_\tau(u,v)=(\tau^2u,\tau v).
\end{displaymath}
Functions with the property
\begin{displaymath}
(f\circ \delta_\tau)(g)=\tau^kf(g)
\end{displaymath}
will be called {\em $\delta_\tau$-homogeneous functions of degree
$k$\/}. The
class of such functions having continuous restriction to the
nonisotropic unit
sphere
 $\Omega^{2n}:=\{(u,v)\in \Space{ H}n|\ u^4+|v|^2=1\}$ is denoted by
$H^k_
\delta(C(\Omega^{2n}))$.

The left and right Haar measure\footnote{The left (right) Haar measure on a
group is a measure that is invariant under the left (right) action of group.}
on the Heisenberg group coincides with
the Lebesgue
measure. The operators of right, left, and two-sided convolution on the
Heisenberg group with kernel $k_{l,r}(g)$ or $k(g_1,g_2)$ are introduced as the
integrals of
the
shift operators $\pi_l(g)$ and $\pi_r(g)$ giving rise to the regular
representation of the Heisenberg group $\Space{ H}n$ on the space
$L^2(\Space{H}n)$:
\begin{eqnarray}
K_{l,r}&=&(2\pi)^{-N/2}\int_{\Space{ H}n}
k_{l,r}(g)\pi_{l,r}(g)\,dg,\label{eq:lrconv}\\
K&=&(2\pi)^{-N}\int_{\Space{ H}n} \int_{\Space{ H}n}
k(g_1,g_2)\pi_l(g_1)\pi_r(g_2)
\,dg_1\,dg_2.\label{eq:tconv}
\end{eqnarray}
where $N=2n+1$.

The Heisenberg group is the simplest non-commutative nilpotent Lie group.
It is well known \cite{MTaylor86}, convolution operators on a
step 2 nilpotent Lie group with kernel $k(g)$ are
pseudodifferential operators (PDO, see
\cite{Hormander85,Shubin78,MTaylor81}) having the following
form:
\begin{eqnarray}\displaystyle
 a(h,D) u(h)& =&(2\pi)^{-N/2} \int_{\Space{R}{N}} e^{i<h,\nu>}
a(h,\nu) \widehat{u}(\nu) d\nu=\nonumber \\ &
=&(2\pi)^{-N} \int_{\Space{R}{N}}  \int_{\Space{R}{N}}
e^{i<h-g,\nu>}\, a(h,\nu)\, u(g)\, dg\,d\nu,
\label{eq:pdo}
\end{eqnarray}
where
\begin{equation}
a(h,\nu)=\widehat{k}( \widetilde{L}_{h} (\nu) ),
\label{eq:symconv}
\end{equation}
and $\widetilde{L}_{h} (\cdot)=\mbox{} ^{t}{L}^{-1}_{h} (\cdot)$ is the
linear operator, which is
inverse and transpose to the operator $L_{h} (\cdot)= -I-
\frac{1}{2}[h,\cdot]$.

In \cite{Shubin78} the more general PDOs containing $\tau$-symbol were
defined:
\begin{equation}
a_{\tau}(h,D) u(h) =(2\pi)^{-N} \int_{\Space{R}{N}}
\int_{\Space{R}{N}}
e^{i<h-g,\nu>}\, a(\tau h+(1-\tau)g,\nu)\, u(g)\, dg\,d\nu.
\label{eq:taupdo}
\end{equation}

If $\tau=1$, then this formula gives the same result as
(\ref{eq:pdo}) and such operator is called PDO with the right symbol. If
$\tau=0$, then this operator is called PDO with the left symbol, if
$\tau=\frac{1}{2}$ --- PDO with the Weyl (symmetric) symbol. The connection
between the
$\tau$-calculus of PDO and the problem of a quantization (in the usual
sense) was  discussed in \cite{Berezin86}.

It is easy to calculate by the formula (\ref{eq:symconv}) not only the
right
symbol of convolution but also any $\tau$-symbol. Indeed  the obvious
equalities (here $(\mbox{\bf ad}h)g= [h,g]$):
\begin{displaymath}
(\mbox{\bf ad} h)h=0, \ (\mbox{\bf ad} h) g=-(\mbox{\bf ad} g)h, \
\end{displaymath}
and (\ref{eq:symconv}) imply:
\begin{displaymath}
L_{h} (h-g) = L_{\tau h+(1-\tau)g}(h-g).
\end{displaymath}
Substituting this to (\ref{eq:taupdo}) one can obtain:
\begin{displaymath}
K u(h) =(2\pi)^{-N} \int_{\Space{R}{N}}
\int_{\Space{R}{N}}
e^{i<h-g,\nu>}\, \widehat{k}(\mbox{} ^{t}{L}^{-1}_{\tau h+(1-\tau)g} \nu)\,
u(g)\, dg\,d\nu.
\end{displaymath}
Thus we have

\begin{prop}\label{pr:tau}\mbox{\cite{Kisil92}}
A $\tau-$symbol of a PDO corresponding to a convolution on a step 2
nilpotent
Lie group does not depend on $\tau$ and has the form:
\begin{equation}
a_{\tau}(h,\nu)=\widehat{k}(\mbox{}^{t}L^{-1}_{h}(\nu)\,)
\label{eq:tausym}
\end{equation}
\end{prop}
\begin{rem}{}\label{re:quantum1}
In the language of quantum mechanics this result
means that the description of a physical system, which symmetry group is
a step 2 nilpotent Lie group, does not depend on a method of the
quantization (right, left or Weyl-symmetric) we use \cite{Berezin86}.
Nilpotent Lie group has special meaning in quantum mechanics, in
particular, the Lie algebra of the Heisenberg group realizes the famous
Heisenberg commutator relations for coordinates and impulses.
\end{rem}

For our purpose we need all irreducible representations of the group
\Heisen{n}. They are given by the Stone-von Neumann theorem
\cite{Kirillov76,MTaylor84,MTaylor86} up to unitary equivalence. For any
$\lambda \in (0,\infty)$ the irreducible noncommutative unitary
representations on $\FSpace{L}{2}(\Space{R}{n})$ are given by
\begin{equation}\label{eq:stone1}
\pi_{\pm\lambda}(t,x,y)=e^{i(\pm\lambda tI\pm\lambda^{1/2}
yM+\lambda^{1/2}xD)},
\end{equation}
where $yM$ and $xD$ are such operators on \FSpace{L}{2}(\Space{R}{n}):
\begin{eqnarray*}\displaystyle
(yM)u(v)&=&\sum y_{j}v_{j}u(v),\\
(xD)u(v)&=&(\frac{1}{i}) \sum x_{j}\frac{\partial u}{\partial v_{j}}.
\end{eqnarray*}
For $(q,p)\in\Space{R}{2n}$, there are also commutative one-dimensional
representations on \Space{C}{}:
\begin{equation}\label{eq:stone2}
\pi_{(q,p)}(t,x,y) u=e^{i(qx+py)}u,\ u\in\Space{C}{}.
\end{equation}
Then relative to \eqref{eq:stone1} -- \eqref{eq:stone2} representations of
convolution algebra are expressed by formulas \cite{MTaylor86}:
\begin{eqnarray}
\pi_{\pm\lambda}[k(t,x,y)]&=&\widehat{k}(\pm\lambda,
\pm\lambda^{1/2}X,\lambda^{1/2}D),\label{eq:cstone1}\\
\pi_{(q,p)}[k(t,x,y)] &=&\widehat{k}(0,q,p)\label{eq:cstone2}.
\end{eqnarray}
The right side of~\eqref{eq:cstone1} specifies a PDO with the Weyl
symbol $\widehat{k}(\pm\lambda,
\pm\lambda^{1/2}x,\lambda^{1/2}\xi)$ accordingly to~\eqref{eq:taupdo} with
$\tau=\displaystyle\frac{1}{2}$. In the right side
of~\eqref{eq:cstone2} one can find just a constant from \Space{C}{}.

\section{The Correspondence Principle between Classical and Quantum
Mechanics}

\begin{figure}[t]
\input{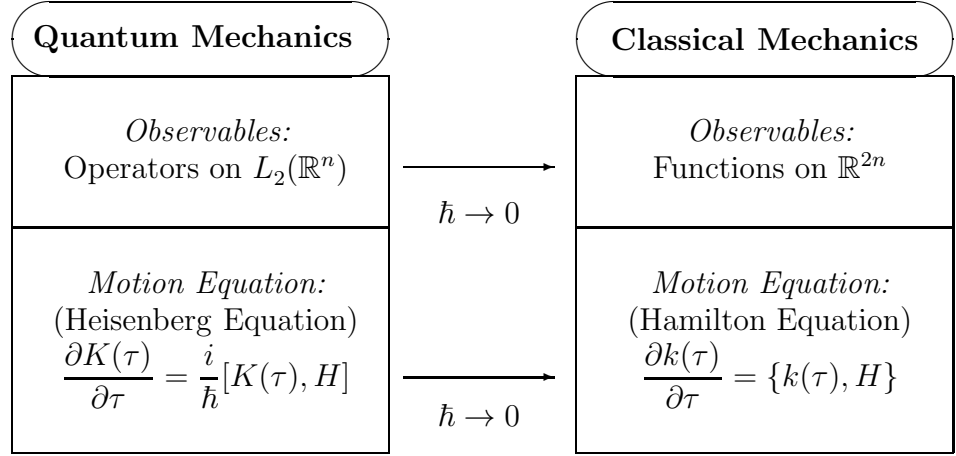}
\caption{The quantization in the usual sense and the correspondence principle
(case of a particle).}\label{fi:usual}
\end{figure}

 The more recent approaches to the quantization problem may be found in
\cite{Berezin74,Berezin75,Unterberger80,Unterberger84}. Let us remind
the sketch of this scheme for the future references (see
Figure~\ref{fi:usual}).

They say that there is a quantization, if
\begin{enumerate}
\item There is a family of operator algebras $\{\algebra{Q}_{\hbar} \,
\mid \hbar\in\Space{R}{}_+\cup \{0\}\}$, where
\begin{enumerate}
\item algebras $\algebra{Q}_\hbar$ for $\hbar\neq 0$ are non-commutative
algebras of operators on some Hilbert spaces;
\item the algebra $\algebra{Q}_\hbar$ for $\hbar =0$ is a commutative
algebra of functions on \Space{R}{2n}.
\end{enumerate}
\item There is a topology on $\Space{R}{}_+$ such that:
\begin{enumerate}
\item there are limits\footnote{\label{fn:exact}The more exact meaning of these
limits
is following: $\displaystyle \lim_{\epsilon\rightarrow 0} a_\epsilon =b$
iff $b$ belong to the closure of the set $\displaystyle
\cup_{t\leq\epsilon}a_t$ for all $\epsilon $.} $\displaystyle
\lim_{\hbar\rightarrow 0}A_\hbar=A_0$, for
$A_\hbar\in\algebra{Q}_\hbar$;
\item for any $A_\hbar$ and $B_\hbar$ the following equalities hold:
\begin{eqnarray}
\lim_{\hbar\rightarrow 0}A_\hbar\circ_\hbar B_\hbar &=& A_0\cdot B_0
\label{eq:hcompos}\\
\lim_{\hbar\rightarrow 0}\frac{i}{\hbar}[A_\hbar,B_\hbar]_\hbar &=&
\{A_0, B_0\}.\label{eq:hcommut}
\end{eqnarray}
Here $\circ_\hbar$ and $\cdot$ in \eqref{eq:hcompos} denote the operator
of composition in $\algebra{Q}_\hbar$ and the ordinary product in
\FSpace{L}{2}(\Space{R}{2n}) correspondingly. In the regular way we denote in
\eqref{eq:hcommut} the commutator of two operators by
$[\cdot,\cdot]_\hbar$, and $\{\cdot,\cdot\}$ is the Poisson brackets of
two functions.
\end{enumerate}
\end{enumerate}

Of course, the accordance between the operator commutator and the
Poisson brackets in~\eqref{eq:hcommut} produce the agreement between the
Heisenberg and Hamilton motion equations when $\hbar\rightarrow 0$ (see
Figure~\ref{fi:usual}).

So we have two very different description for one system with only a
thin bridge among them: the limit by $\hbar\rightarrow 0$. Please, do
not ask the question:{\em Why should we consider such limits by
$\hbar$, if $\hbar$ is the Planck {\bfseries constant\/}?\/} I think,
such
questions cannot be answered within mentioned scheme.

\begin{figure}[p]
\input{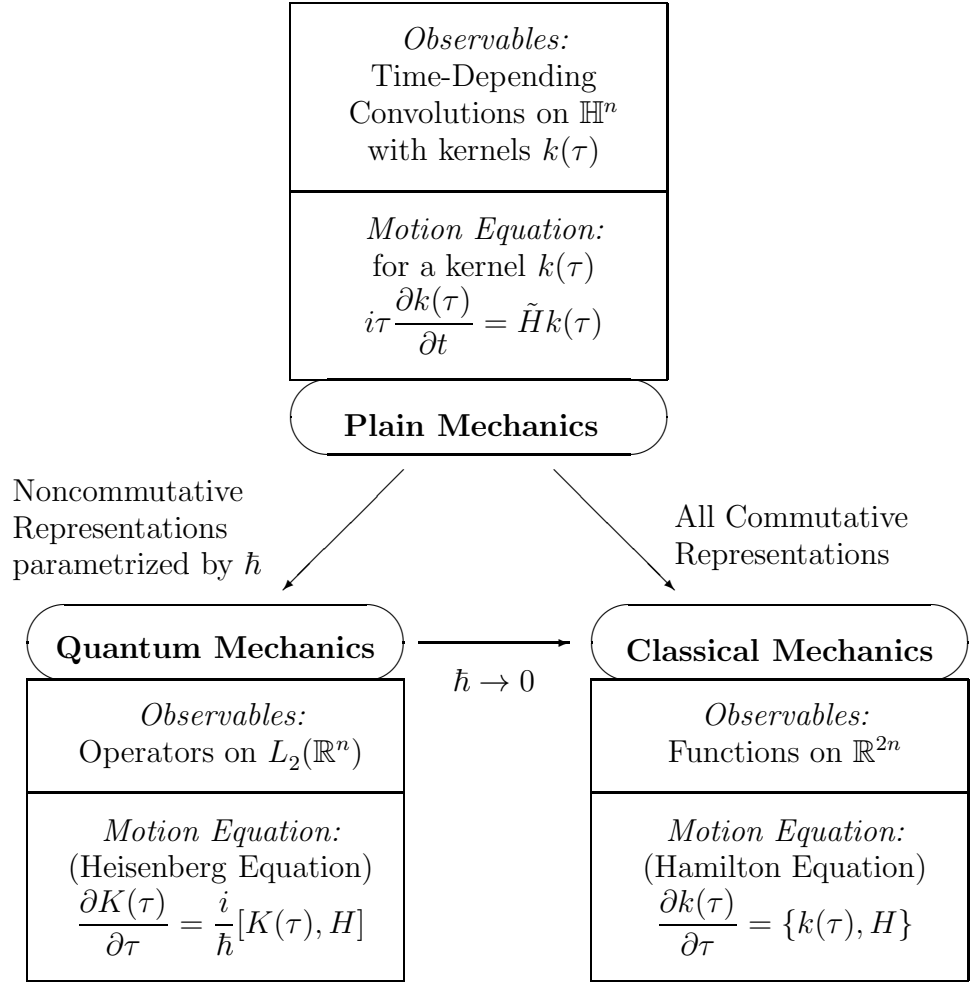}
\caption{Plain Mechanics: Superstructure for the usual scheme (case of a
particle).}\label{fi:plain}
\end{figure}

\section{Joining of Classical and Quantum Mechanics}\label{se:plain}

Now we can discuss a new approach to our question.

We will speak
that there is a {\em plain mechanical description\/} of a system if the
following
requirements hold:
\begin{enumerate}
\item There is an operator algebra \algebra{P} which has
\begin{enumerate}
\item the family of all infinite-dimensional noncommutative irreducible
representations parametrized by points of a set $P$:
\begin{equation}
\pi_\hbar: \algebra{P} \rightarrow \algebra{P}_\hbar,\ \hbar\in P,\
\object{dim}(\algebra{P}_{\hbar})=\infty.
\end{equation}
 We will call $P$ by the {\em set of the Planck constants\/};
\item the family of all (commutative) one-dimensional representations
parametrized by points of a set $M$. This family gives the mapping of
\algebra{P} into algebra of functions on $M$ by the obvious rule:
\begin{equation}\label{eq:pizero}
\pi_x: P\rightarrow p(x)\in\Space{C}{},\ P\in\algebra{P},\ x\in M
\end{equation}
The set $M$ will be called as {\em the phase space of the classical  system},
and the mapping defined by~\eqref{eq:pizero} from \algebra{P} to an algebra of
functions on $M$ will be denoted by $\pi_0$;
\item the topology\footnote{This topology should be naturally generated
by the structure of the algebra \algebra{P}, for example it may be the Jacobson
topology~\cite{Dixmier69}  or the *-bundle
topology~\cite{DaunHof68,Hofmann72}.} $T$ on the set $P\cup M$ of
all  its representations.
\end{enumerate}

\item\label{it:commute} The algebra \algebra{P} is equipped by the operation
$[\cdot,\cdot]$ of the commutation such that its image under mapping
$\pi_0$ should coincide with the Poisson brackets\footnote{The Poisson
brackets assumes that $M$ should have the structure of a manyfold. We
don't discuss now how it is happened in the general case. Luckily, this
structure will appear very natural in Subsection~\ref{ss:qplain}.}  on
$M$.

\item The diagram on Figure~\ref{fi:plain} should be commutative. The left
down-going arrow denotes the set of  all infinite-dimensional representations
of
\algebra{P}, the right arrow indicates the mapping $\pi_0$ and the
horizontal arrow means the limit in topology $T$.
\end{enumerate}

Comparing the Figure~\ref{fi:usual} and Figure~\ref{fi:plain} it is easy to
see that plain mechanics is a superstructure on the usual scheme of a
quantization.

\begin{rem}\label{re:measurement} To give a short philosophical
interpretation of plain mechanics I would like to stress the following.
By my opinion, plain mechanics correspond to the inner structure of the
world. But we cannot see this inner structure. During the process of an
observation (measurement) a representation of the world is selected
from all possible ones. It may be either the classic representation or any from
different (for different $\hbar$)  quantum ones. Which representation
was selected  depends on the observer and his equipment (apparatus).
\end{rem}

\section{Example: a Particle in $n$-dimensional Space}
In this section we give an illustration of the given abstract schemes
on the simplest example of a particle in $n$-dimensional space. The Weyl
and the Wick--Berezin quantization (Subsections~\ref{ss:qWeyl}
and~\ref{ss:qWick}) represent the usual
methodology. In Subsection~\ref{ss:qplain} a realization of the abstract
scheme from Section~\ref{se:plain} by two-sided convolutions on
the Heisenberg group is done.

There is a simple reminding of the classical case for our problem. States
$(q,p)$ of a particle in \Space{R}{n} form the manifold (phase
space) \Space{R}{2n}. The observables are real functions on  the phase
space \Space{R}{2n} and the value of measurement an observable $k$ in a
state $(q,p)\in\Space{R}{2n}$ is just $k(q,p)$. The change of an observable
$k(\tau)$ during the time\footnote{We employ the symbol $\tau$ for notation
of time-parameter because the letter $t$ have been already used for the first
coordinate on the Heisenberg group.} defined by the Hamilton equation:
\begin{equation}\label{eq:Hamilton}
\frac{\partial k(\tau)}{\partial \tau}=\{k(\tau),H\},
\end{equation}
where $H(q,p)$ is the Hamilton function of the full energy of the particle.

\subsection{The Weyl Quantization: PDO Calculus on \Space{R}{n}
}\label{ss:qWeyl}

 The realization is based on the well-known symbolic
calculus of PDO \cite{Hormander85,Shubin78,MTaylor81}. States are the
functions $f(x)$ from $\FSpace{L}{2}(\Space{R}{n})$ and observables are
operators on $\FSpace{L}{2}(\Space{R}{n})$. The mathematical expectation
of an observable $K$ on a state $f(x)$ is equal to
$\scalar{Kf}{f}$\footnote{This means the scalar product on
$\FSpace{L}{2}(\Space{R}{n})$.}. The Heisenberg equation
\begin{equation}\label{eq:Heisenberg}
\frac{\partial K(\tau)}{\partial \tau}=\frac{i}{\hbar}[H,K(\tau)]
\end{equation}
defines the motion of observables. The classic observables of a
coordinate $q_i$ and an impulse $p_i$ correspond to the operators of
multiplication by $x_i$ and derivative $\displaystyle
\frac{\hbar}{i}\frac{\partial
}{\partial x_i}$ correspondingly.

The plainness of this description is suspended by the question:
{\itshape Which operator $(x_i \displaystyle \frac{\partial }{\partial
x_i},\ \displaystyle \frac{\partial }{\partial x_i} x_i \hbox{ or }
\displaystyle \frac{1}{2}(x_i\frac{\partial }{\partial x_i} +
\frac{\partial }{\partial x_i}x_i))$ corresponds to the classical
observable $p_iq_i$?\/} Different answers on this question generate
different correspondences between symbols and operators in the calculus
of PDO (see \eqref{eq:taupdo}). Remark~\ref{re:quantum1} defines an
important class of symbols without this obstacle.

\subsection{The Berezin (Anti-Wick) Quantization: Toeplitz Operators in the
Fock Space}\label{ss:qWick}

Again we give only short summary of this topic, the relevant information
may be found in
\cite{Berezin74,Berezin75,Berezin86,Coburn90,Guillemin84}. Let
$\FSpace{L}{2}(\Space{C}{n},d\mu_{n})$ be a space of all square-integrable
functions on
\Space{C}{n} with respect to the Gaussian measure
\begin{displaymath}
d\mu_{n}(z)=\pi^{-n}e^{-z\cdot\overline{z}}dv(z),
\end{displaymath}
where $dv(z)=dxdy$ is the usual Euclidean volume measure on
$\Space{C}{n}=\Space{R}{2n}$. Denote
by $P_{n}$ the orthogonal Bargmann projector \cite{Bargmann61}  of
$\FSpace{L}{2}(\Space{C}{n},d\mu_{n})$ onto the Fock space
$\FSpace{F}{2}(\Space{C}{n})$, the subspace of
$\FSpace{L}{2}(\Space{C}{n},d\mu_{n})$ consisting of all entire
functions. Then the formula
\begin{equation}
k(q,p)\rightarrow T_{k(q+ip)}=P_n k(p+iq)I
\end{equation}
defines an another (anti-Wick or Berezin) quantization, which maps a function
$k(q,p)$ on \Space{R}{2n} to the Toeplitz operator $T_k$  with the
pre-symbol $k(q+ip)$ on \Space{C}{n}. There is an identification
between the
Berezin quantization and the Weyl quantization
\cite{Berezin74,Coburn90,Guillemin84}.

\subsection{Plain Mechanics: Two-Sided Convolutions on the
Heisenberg Group}\label{ss:qplain}

Now we illustrate the scheme from Section~\ref{se:plain}. Let us take a
convolution operator algebra on the Heisenberg group. The kernels of
convolutions may be taken from $\FSpace{L}{1}(\Heisen{n})$ for
example. Using the formulas \eqref{eq:cstone1} --\eqref{eq:cstone2} we
conclude that the set of the Planck constants $P$  coincides with $\Space{R}{}
\setminus 0$.
The phase space in the sense of Section~\ref{se:plain} (the set of all
one-dimensional representations) complies with the classical phase space
\Space{R}{2n} and we can transfer the manifold structure of
\Space{R}{2n} to $M$ (with the associated Poisson brackets). The Jacobson
topology on $P$ is induced by the usual topology of the real line and
any interval $(-\alpha, 0) \hbox{ or } (0, \alpha)\subset P,\ \alpha >0$ is
everywhere dense in $M$. This mean
that the limit
\begin{equation}\label{eq:limit}
\widehat{k}(\pm\lambda,\pm\lambda^{1/2}X,\lambda^{1/2}D)\rightarrow
\widehat{k}(0,q,p)
\end{equation}
while $\lambda\rightarrow 0$ is well defined\footnote{See
footnote~\ref{fn:exact} for the exact definition of this limit.} in the
Jacobson topology.

Now we check the commutator property (see item~\ref{it:commute} from
page~\pageref{it:commute}):
\begin{prop}\label{pr:commutator}
The limit of $\displaystyle \frac{i}{\lambda}\pi_\lambda([K_1,K_2])$ by
$\lambda\rightarrow 0$ is
equal to $\{\widehat{k}_1(0,q,p),\widehat{k}_2(0,q,p)\}$. Here
$\pi_\lambda$ is defined by~\mbox{\rm\eqref{eq:cstone1}}, $K_1$ and $K_2$ are
convolutions with kernels $k_1, k_2$ respectively.
\end{prop}
\begin{proof}
Using the standard PDO calculus one can calculate that the symbol of
commutator image under representation $\pi_\lambda$ is equal to
\begin{eqnarray*}
\object{sym}(\pi_\lambda([K_1,K_2]))(q,p)&=&-i\frac{\partial
\widehat{k}_1(\lambda,\lambda^{1/2}q,\lambda^{1/2}p)}{\partial
q}\frac{\partial
\widehat{k}_2(\lambda,\lambda^{1/2}q,\lambda^{1/2}p)}{\partial p}\\
&&+i
\frac{\partial
\widehat{k}_1(\lambda,\lambda^{1/2}q,\lambda^{1/2}p)}{\partial
p}\frac{\partial
\widehat{k}_2(\lambda,\lambda^{1/2}q,\lambda^{1/2}p)}{\partial q}\\
&&+\hbox{ (derivatives of orders $> 2$)}
\end{eqnarray*}
Note that in this expansion any derivative of order $m$ has the
vanishing order $\lambda^{m/2}$ when $\lambda\rightarrow 0$. Thus if we
multiply the image of commutator $\object{sym}(\pi_\lambda([K_1,K_2]))$ by
$\displaystyle \frac{i}{\lambda}$ and
take the limit accordingly to~\eqref{eq:limit} then we obtain the
assertion.
\end{proof}
\begin{rem}
It is easy to see that many essentially different observables may have
the same classical representation. Really, if
\begin{equation}
\widehat{k}_1(0,q,p)=\widehat{k}_2(0,q,p)
\end{equation}
then observables the $K_1$ and $K_2$ are identical from the classical point
of view. A differentness among them can be found only on a quantum level
after selection a tool for the observation (see
Remark~\ref{re:measurement}). This note gives another dimension to the
old dispute on the existence of hidden variables.
\end{rem}
\begin{rem}
If we take the algebra of observables with $\delta_\tau$-homogeneous
kernels  $\widehat{k}(t,x,y)\in H_\delta^0(C(\Omega^{2n}))$ (see
Section~\ref{se:background}) then different quantum representations will
not depend from $\hbar$. Nevertheless the classical limit of these ``constant''
quantum observables by $\hbar\rightarrow 0$
is still  defined by~\eqref{eq:limit}.
\end{rem}

Let $H$ be a convolution on the Heisenberg group corresponding to the
Hamiltonians $H_\lambda$ by representations $\pi_\lambda: H \rightarrow
H_\lambda$. Let us calculate the kernel  $c(h)$ of its commutator $[H,K]$ with
a convolution $K$. We denote the

ernels of operators $H$ and $K$ by $j(h)$ and $k(h)$ correspondingly. It is
easy to verify that:
\begin{eqnarray*}
c(h)&=&\int_{\Heisen{n}}(j(g)k(h*g)- k(g)j(h*g))\,dg\\
    &=&\int_{\Heisen{n}}j(g)k(h*g)\,dg - \int_{\Heisen{n}} k(g)j(h*g)\,dg \\
    &=&\int_{\Heisen{n}}j(g)k(h*g)\,dg - \int_{\Heisen{n}} k(g)j(h*g)\,dg \\
    &=&\int_{\Heisen{n}}j(g)k(h*g)\,dg - \int_{\Heisen{n}}
k(g_{1}^{-1}*h)j(h*g_{1}^{-1}*h)\,dg_{1} \\
    &=&\int_{\Heisen{n}}\int_{{\Heisen{n}}}j(g)\delta (g_{1})
k(g_{1}^{-1}*h*g)\,dgdg_{1}\\
    & &- \int_{\Heisen{n}}\int_{\Heisen{n}}\delta(g)j(h*g_{1}^{-1}*h)
k(g_{1}^{-1}*h*g)\,dgdg_{1} \\
    &=& \int_{\Heisen{n}}\int_{\Heisen{n}}(j(g)\delta (g_{1})
-\delta(g)j(h*g_{1}^{-1}*h))   k(g_{1}^{-1}*h*g)\,dgdg_{1}. \\
\end{eqnarray*}
Here $\delta (g)$ is the Dirac function. Thus
\begin{equation}\label{eq:commutator}
[H,K] = \widetilde{H}K,
\end{equation}
where $\widetilde{H}$ is an operator of two-sided convolution with the kernel
depended from a point of $h\in\Heisen{n}$:
\begin{equation}\label{eq:hamkern}
\widetilde{j}(h,g,g_1)=(j(g)\delta (g_{1}) -\delta(g)j(h*g_{1}^{-1}*h)).
\end{equation}
Such operator $\widetilde{H}$ belongs to an algebra generated by two-sided
convolutions~\eqref{eq:tconv}  and operators of multiplication by functions on
\Heisen{n}~\cite{Kisil92}.
\begin{prop}\mbox{\itshape\bfseries(The Main Equation of Plain Mechanics)}
The equation
\begin{equation}\label{eq:main}
i\tau\frac{\partial K(\tau)}{\partial t}=-\widetilde{H}K(\tau),
\end{equation}
where $\widetilde{H}$ is an operator defined by the
kernel~{\/\mbox{\rm\eqref{eq:hamkern}}} turns into the Heisenberg equation
under a mapping $\pi_\lambda$ and turns into the Hamilton equation under the
mapping $\pi_0$.
\end{prop}
\begin{proof} Let convolution $K(\tau)$ have the kernel $k(\tau,t,x,y)$. Then
convolution $i\displaystyle \tau\frac{\partial K(\tau)}{\partial t}$ has the
kernel $i\displaystyle \tau\frac{\partial k}{\partial t}(\tau,t,x,y)$. Under
representations \eqref{eq:cstone1} --\eqref{eq:cstone2} the convolution $i\displaystyle
\tau\frac{\partial K(\tau)}{\partial t}$ has the images
\begin{eqnarray*}
i\lambda\frac{\partial \widehat{k}}{\partial
\tau}(\tau,\pm\lambda,\pm\lambda^{1/2}X,\lambda^{1/2}D),&&
i\frac{\partial \widehat{k}}{\partial \tau}(\tau,0,q,p)
\end{eqnarray*}
correspondingly. We take a liberty to denote by $\tau$  both the time-parameter
and  its dual in the Fourier transform sense. Note please that $\lambda$ have
the meaning of the Planck constant here.

The right side of~\eqref{eq:main} in accordance with~\eqref{eq:commutator} is
equal to commutator $[H,K]$. Thus taking in account
Proposition~\ref{pr:commutator} we have
\begin{eqnarray}
\pi_\lambda&:&i\tau\frac{\partial K(\tau)}{\partial
t}=\widetilde{H}K(\tau)\nonumber \\
&&\rightarrow i\lambda\frac{\partial \widehat{k}}{\partial
\tau}(\tau,\pm\lambda,\pm\lambda^{1/2}X,\lambda^{1/2}D)=[H_\lambda,
\widehat{k}(\tau,\pm\lambda,\pm\lambda^{1/2}X,
\lambda^{1/2}D)]\label{eq:iheisenberg}\\
\pi_0&:&i\tau\frac{\partial K(\tau)}{\partial t}=\widetilde{H}K(\tau)\nonumber
\\
&&\rightarrow i\frac{\partial \widehat{k}}{\partial \tau}(\tau,0,q,p)=
i\left( \frac{\partial \widehat{k}}{\partial q}\frac{\partial
\widehat{j}}{\partial p} -\frac{\partial \widehat{j}}{\partial q}
\frac{\partial \widehat{k}}{\partial p}\right)(\tau,0,q,p).\label{eq:ihamilton}
\end{eqnarray}
Here $H_\lambda$ is the image of the convolution $H$ under the representation
$\pi_\lambda$ and $j$ is the kernel of the convolution $H$. It is clear that
equations~\eqref{eq:iheisenberg} -- \eqref{eq:ihamilton} coincide with the
Heisenberg equation~\eqref{eq:Heisenberg} and the Hamilton
equation~\eqref{eq:Hamilton} correspondingly. 
\end{proof}
\begin{rem}
The Heisenberg equation (equation for observables)~\eqref{eq:main} looks like
the Schr\"odinger equation (equation for states) in quantum mechanics. The left
side of~\eqref{eq:main} contains the partial derivative by $t$, so it looks
like $t$ is the time

parameter in plain mechanics.
\end{rem}

\bibliographystyle{amsplain}
\bibliography{abbrevmr,akisil,aphysics,analyse}

\end{document}